\documentclass[a4paper,10pt]{article}
\usepackage{thesis}
\usepackage{natbib}
\usepackage{graphicx}
\usepackage{transparent}
\usepackage[font={singlespacing,small,it}]{caption}
\begin{document}
\font\myfont=cmr12 at 15pt
\title{{\myfont KdV cnoidal waves in a traffic flow model with periodic boundaries}}
\date{}
\author{Laura Hattam \footnote{University of Reading, l.hattam@reading.ac.uk}}
\maketitle

\abstract{An optimal-velocity (OV) model describes car motion on a single lane road. In particular, near to the boundary signifying the onset of traffic jams, this model reduces to a perturbed Korteweg-de Vries (KdV) equation using asymptotic analysis. Previously, the KdV soliton solution has then been found and compared to numerical results (see \citet{mur99}). Here, we instead apply modulation theory to this perturbed KdV equation to obtain at leading order, the modulated cnoidal wave solution. At the next order, the Whitham equations are derived, which have been modified due to the equation perturbation terms. Next, from this modulation system, a family of spatially periodic cnoidal waves are identified that characterise vehicle headway distance. Then, for this set of solutions, we establish the relationship between the wave speed and the modulation term, which is dependent upon the number of oscillations over the spatial domain. This analysis is confirmed with comparisons to numerical solutions of the OV model.}
\section{Introduction}
The study of traffic flow has uncovered some interesting phenomena such as the propagation of nonlinear density waves representing congestion. To determine traffic behaviour, a variety of modelling techniques are used, which include the application of car following, cellular automation, gas kinetic and hydro dynamical models. Refer to \citet{nag02} for a discussion of the different methods.

Here, we concentrate on a car following model that governs vehicle motion on a single lane road with periodic boundaries. The OV model proposed by \citet{new61} is applied, which is
\beq
\frac{d x_j(t+\tau)}{dt}=V(x_{j+1}(t)-x_j(t)),\label{newell}\eeq
where $x_j(t)$ is the position of car $j$ at time $t$, $\tau$ is the delay time of the driver, $V$ is the optimal velocity and $j=0,1,2,\ldots,N$ for $N$ cars on the road.

We will examine (\ref{newell}) when it reduces to a perturbed KdV equation. This occurs only within a certain stability zone, which is referred to as `metastable'. It is well-known that the unperturbed KdV equation has the travelling wave solution, the cnoidal wave of modulus $m$, where $m\in(0,1)$ (see \citet{kdv95}). For $m\rightarrow 1$, this becomes the soliton solution.

\citet{mur99} explored this traffic flow problem with open boundaries. They derived a perturbed KdV equation from an OV model and then obtained the KdV soliton solution. This result was compared to numerical simulations with good agreement. As well, these solutions were shown to disappear after some time. Later, \citet{zhu08} performed a similar analysis for periodic boundaries. Numerically they found large amplitude waves with narrow peaks of both upward and downward form, which they referred to as KdV solitons. Nonlinear analysis was also applied to find the soliton solution. The focus of these investigations into the metastable dynamics was the soliton.

Additionally, \citet{yu10} and \citet{zho14} performed numerical examinations of traffic OV models, where similar nonlinear behaviour emerged. However, the simulations also revealed steady travelling waves in the metastable region that were of a similar form to the KdV modulated cnoidal wave. Therefore, other KdV solutions besides the soliton can occur in OV models. Hence, further asymptotic work is needed to establish the connection between the numerically observed travelling density waves and the nonlinear analysis. To achieve this, we will apply modulation theory to the perturbed KdV equation so that the entire family of possible solutions to the reduced traffic model are obtained.

\citet{whi74} developed modulation theory for the KdV equation, which was a multi-scale technique that gave modulated wavetrain solutions. A system of first order partial differential equations describing the modulations were also found, now known as the `Whitham Equations'. \citet{gur87} extended the modulation theory to include Burgers damping, where the Whitham equations with additional terms to account for friction were formed. \citet{myi95} instead incorporated an arbitrary damping term, and then considered three different types of damping to analyse the subsequent wavetrain solutions. The monograph by \citet{kam00} consisted of a comprehensive derivation of the KdV Whitham Equations. These previous studies will be followed here to obtain the modified Whitham equations, which correspond to our traffic problem.

This paper concentrates on the identification of steady travelling wave solutions of the modulation equations since these are found numerically. More specifically, in Section $2$, we outline the traffic OV model and the asymptotic analysis used to then obtain a perturbed KdV equation. Next, in Section $3$, a multi-scale perturbation technique is applied to this equation and as a result, the modulation equations are derived. Then, steady solutions of these equations are pursued in Section $4$. This asymptotic theory is next related to the traffic problem in Section $5$. Lastly, in Section $6$, numerical simulations of the OV model are depicted and compared to the asymptotic solutions.

\section{Traffic Flow Model}
We outline the transformation of (\ref{newell}) into a perturbed KdV equation within the metastable zone. Firstly, undertaking a Taylor series expansion of (\ref{newell}) and expressing this in terms of the vehicle headway, $\Delta x_j=x_{j+1}-x_j$, gives
%
\beq
\frac{d^2 \Delta x_j}{d t^2}=\hat{a}\left(V(\Delta x_{j+1}(t))-V(\Delta x_{j}(t))-\frac{d \Delta x_j}{d t}\right),\label{traf_hw}
\eeq
where $\hat{a}=1/\tau$ is the drivers sensitivity. This model was proposed by \citet{ban95}.

As well, the following optimal velocity function defined by \citet{ban95} is used
\beq
V(\Delta x_j(t))=\frac{v_{max}}{2}\left(\tanh(\Delta x_j-h_c)+\tanh(h_c)\right),\label{veldef}\eeq
where $h_c$ is the safety distance and $v_{max}$ is the maximal velocity. We choose $v_{max}=2$ and $h_c=4$ for convenience.

\citet{ge05} outlined three OV models that describe car motion. Here, their Model B is applied, which uses (\ref{traf_hw}) and (\ref{veldef}), and has the linear stability criteria
\beq
\tau\le\tau_s=\frac{1}{2V^{'}(h)}.\label{linstab}\eeq
If this condition is satisfied, then the steady solution to (\ref{traf_hw}) is stable. This steady state is $\Delta x_j(t)=h$, where $h$ is the uniform headway. The curve given by $\tau=\tau_s$ is labelled the `neutral stability line' as it represents the boundary between no traffic jams and jams.

Moreover, \citet{ge05} detailed the application of an asymptotic method to reduce (\ref{traf_hw}) to a perturbed KdV equation, using the change of variables
\beq
x=-\epsilon\sqrt{\frac{6}{V^{'}(h)}}(j+V^{'}(h)t),\quad \bar{t}=\epsilon^3\sqrt{\frac{6}{V^{'}(h)}} t,\quad \epsilon^2=1-(\tau/\tau_s),\quad 0<\epsilon\ll 1,\label{jttoxtb}\eeq
and letting
\beq
\Delta x_j(t)=h+\frac{\epsilon^2}{V^{''}(h)} u.\eeq

Consequently, (\ref{traf_hw}) becomes
\begin{equation}
u_{\bar{t}}+\nu uu_x+\lambda u_{xxx}+\epsilon \bar{V}(u)=0,\label{meqnuns}\eeq
where
\beq
\bar{V}(u)=\mu u_{xx}+\gamma u_{xxxx}+\eta\left(u^2\right)_{xx},\label{defVb}
\eeq
and
\beq
\lambda=1,\quad\nu=1,\quad \mu=-\sqrt{\frac{3V^{'}(h)}{2}},\quad \gamma= \frac{3}{2}\sqrt{\frac{3}{2V^{'}(h)}},\quad \eta=\frac{1}{2}\sqrt{\frac{3}{2V^{'}(h)}}.\label{par_val}\eeq
This is a perturbed KdV equation as $\epsilon$ is small. We have introduced here the parameters $\lambda,\;\nu,\;\mu,\;\gamma,\;\eta$ so that the perturbation analysis in Sections $3$ and $4$ is generalised. Then, this is related to our traffic problem in Sections $5$ and $6$.

Since $0<\epsilon\ll 1$, $\tau$ is chosen such that the solutions to (\ref{meqnuns}) are positioned very close to the neutral stability line, however, they will still satisfy (\ref{linstab}). Therefore, as $t$ becomes very large, the headway will tend to steady state $h$. This region is classified as metastable. Soliton density waves within this zone have been shown numerically to propagate for long times and eventually disappear (refer to \citet{mur99}). We however will demonstrate that a large set of long-time persisting cnoidal wave solutions exist. Here, it is only periodic boundaries that are considered, therefore spatially periodic solutions are sought. It should be noted that reference throughout to the spatial domain corresponds to $j\in[0,N]$. So, to implement periodic boundary conditions, we ensure
\beq
\Delta x_0(t)=\Delta x_N(t),\quad \left.\frac{\partial \Delta x_{j}(t)}{\partial j}\right|_{j=0} = \left.\frac{\partial \Delta x_{j}(t)}{\partial j}\right|_{j=N},\quad t\geq 0.\label{pbcon}\eeq

\section{Perturbation Analysis and the Modulation Equations}
A multi-scale perturbation approach that follows \citet{myi95} is applied to (\ref{meqnuns}). The monograph by \citet{whi74} used this perturbation technique for the analysis of the unperturbed system ((\ref{meqnuns}) with $\epsilon\rightarrow 0$) and \citet{myi95} extended this theory to include an arbitrary damping term.

To begin, let
\bdis
u\left(x,\bar{t}\right)=u_0(\theta,X,T)+\epsilon u_1(\theta,X,T)+\epsilon^2u_2(\theta,X,T)+\ldots,\edis
and introduce the variables
\begin{equation}
\theta=\frac{1}{\epsilon}\Theta(X,T),\quad X=\epsilon x,\quad T=\epsilon \bar{t}.\label{c5pertvar}\eeq
Next, the frequency $\omega$, wave number $k$, and wave speed $c$, are given by
\begin{equation}
\omega=-\Theta_T,\quad k=\Theta_X,\quad\omega=kc.\label{c5pertvar2}\eeq
For $\Theta_{XT}=\Theta_{TX}$, it is required
\begin{equation}
k_T+\omega_X=0.\label{concon}\eeq
Making this change of variables, (\ref{meqnuns}) becomes at first and second order
\bsub
\label{12orde}
\baln
&O(1): -cu_{0,\theta}+\lambda k^2 u_{0,\theta\theta\theta} +\nu u_0u_{0,\theta}=0,\label{12ordea}\\
&O(\epsilon): -k c u_{1,\theta}+\lambda k^3 u_{1,\theta\theta\theta}+\nu k(u_0u_1)_{\theta}+g=0,\label{12ordeb}
\end{align}
\esub
where
\begin{equation}
g=u_{0,T}+\nu u_0u_{0,X}+3\lambda(k^2 u_{0,\theta\theta X}+kk_Xu_{0,\theta\theta})+\bar{V}(u_0).\label{c5hdef}\eeq

The solution to (\ref{12ordea}) is
\begin{equation}
u_0=ab+d+a\operatorname{cn}^2(\beta(\theta-\theta_0);m),\label{u0defc5}\eeq
where
\bsub
\label{abddef}
\begin{align}
a=&\frac{12\lambda}{\nu}(mk\beta)^2,\label{abddefa}\\
d=&\frac{c}{\nu}-\frac{a}{3m^2}\left(2-m^2-3\frac{E(m)}{K(m)}\right),\label{abddefb}\\
b=&\frac{1-m^2}{m^2}-\frac{E(m)}{m^2 K(m)},\label{abddefc}\\
\beta=&K(m)/P.\label{abddefd}
\end{align}\esub
This is the cnoidal wave solution with period $2P$ in $\theta$, where $P$ is a fixed constant. The parameters $m,\;a,\;b,\;d,\;k,\;\theta_0,\;\beta$ are slowly varying, dependent on the slow variables $X$ and $T$. As well, $b$ is chosen so that the mean value of $u_0$ is $d$. The function $\operatorname{cn}$ is the Jacobi elliptic function and $K(m)$, $E(m)$ are the elliptic integrals of the first and second kind respectively.

Integrating (\ref{12ordea}) twice, we find
 \begin{equation}
\frac{\lambda k^2}{\nu}u_{0,\theta}^2=2\hat{D}+2\hat{C}u_0+Uu_0^2-\frac{1}{3}u_0^3,\label{1oi2}\eeq
where $\hat{C}$ and $\hat{D}$ are integration constants and $U=c/\nu$. These constants may be expressed in terms of the cnoidal wave parameters,
\bsub
\label{cdeqn}
\baln
&\hat{C}=-\frac{1}{3}(ab+d)^3+\frac{1}{2}U(ab+d)^2-\frac{a^2}{6m^2}(1-m^2)(ab+d),\label{cdeqna}\\
&\hat{D}=\frac{1}{2}(ab+d)^2-U(ab+d)+\frac{a^2}{6m^2}(1-m^2).\label{cdeqnb}
\end{align}
\esub

To ensure $u_1$ is periodic in $\theta$ over $2P$, \citet{myi95} impose the following integral conditions
\begin{equation}
\int_{-P}^Pg d\theta=0,\quad\int_{-P}^Pu_0 g d\theta=0,\label{c5pcond}\eeq
where $g$ is given by (\ref{c5hdef}). Written in full, (\ref{concon}) and (\ref{c5pcond}) take the form
\bsub
\label{begneqns}
\baln
&\frac{\partial d}{\partial T}+\frac{\partial\;}{\partial X}(\nu(Ud+\hat{C}))+\frac{1}{2P}\int_{-P}^P\bar{V}(u_0)d\theta=0,\label{begneqnsa}\\
&\frac{\partial \;}{\partial T}(Ud+\hat{C})+\frac{\partial\;}{\partial X}(\nu U(Ud+\hat{C})-\nu\hat{D})+\frac{1}{2P}\int_{-P}^Pu_0\bar{V}(u_0)d\theta=0,\label{begneqnsb}\\
&k_T+\nu(kU)_X=0.\label{begneqnsc}
\end{align}
\esub
The manipulation of (\ref{1oi2}) was used to express the integrals in (\ref{c5pcond}), in terms of $\hat{C},\;\hat{D},\;U$ and $d$. Note that $u_0$ and its derivatives with respect to $\theta$ are assumed to be periodic over $\theta\in[-P,P]$. For further detail, see \citet{myi95}.

\citet{kam00} outlined the derivation of the Whitham modulation equations for the KdV equation ((\ref{meqnuns}) with $\epsilon\rightarrow 0$). This is a third order system written in terms of the slowly varying Riemann invariants $r_1$, $r_2$, $r_3$. To transform (\ref{begneqns}) into the modulation equations, these workings were followed except the effects of the perturbation terms, $\bar{V}$, were incorporated into the analysis. After a number of steps, which have been omitted here, (\ref{begneqns}) reduces to the modulation equations
\begin{equation}
\frac{\partial r_i}{\partial T}+Q_i\frac{\partial r_i}{\partial X}=M_i,\label{modsys}\eeq
where
\bdis
Q_i=\nu U+\frac{\nu}{6}\frac{k}{\partial_{r_i}k},\quad i=1,2,3,\quad U=\frac{1}{6}(r_1+r_2+r_3),\edis
and
\bdis\bspl
M_1=&-\frac{k}{2P(\partial_{r_1}k)(r_3-r_1)(r_2-r_1)}\times\\
&\quad\left(\int_{-P}^P u_0\bar{V}(u_0)d\theta -\frac{1}{2}(r_2+r_3-r_1)\int_{-P}^P \bar{V}(u_0)d\theta\right),\\
M_2=&-\frac{k}{2P(\partial_{r_2}k)(r_3-r_2)(r_1-r_2)}\times\\
&\quad\left(\int_{-P}^P u_0\bar{V}(u_0)d\theta -\frac{1}{2}(r_1+r_3-r_2)\int_{-P}^P \bar{V}(u_0)d\theta\right),\\
M_3=&-\frac{k}{2P(\partial_{r_3}k)(r_1-r_3)(r_2-r_3)}\times\\
&\quad\left(\int_{-P}^P u_0\bar{V}(u_0)d\theta -\frac{1}{2}(r_1+r_2-r_3)\int_{-P}^P \bar{V}(u_0)d\theta\right).\end{split}\edis
This system is consistent with \citet{myi95} and decribes the slow variation of the leading order solution (\ref{u0defc5}), which is, written in terms of $r_i$,
\bdis
u_0(\theta,X,T)=\frac{1}{2}(r_1+r_3-r_2)+(r_2-r_1)\operatorname{cn}^2\left(\sqrt{\frac{\nu(r_3-r_1)}{12\lambda k^2}}(\theta-\theta_0);m\right),\edis
where $m^2=(r_2-r_1)/(r_3-r_1)$ and $r_i$ are dependent on $X$ and $T$.

Now, we apply this analysis to the perturbed KdV equation outlined in Section $2$, where $\bar{V}$ is defined by (\ref{defVb}). So,
\bdis\bspl
\frac{1}{2P}\int_{-P}^P u_0\bar{V}(u_0)d\theta=&\frac{1}{2P}\int_{-P}^P u_0(\mu k^2u_{0,\theta\theta}+\gamma k^4 u_{0,\theta\theta\theta\theta}+\eta k^2(u_0^2)_{\theta\theta})d\theta\\
=&-\mu k^2 \frac{1}{2P}\int_{-P}^P u_{0,\theta}^2d\theta+\gamma k^4 \frac{1}{2P}\int_{-P}^P u_0u_{0,\theta\theta\theta\theta}d\theta\\
&+\eta k^2\frac{1}{2P}\int_{-P}^Pu_0\left(-\frac{2\lambda k^2}{\nu}u_{0,\theta\theta\theta\theta}+2Uu_{0,\theta\theta}\right) d\theta\\
=&-\left(\mu k^2+2U\eta k^2\right) \frac{1}{2P}\int_{-P}^P u_{0,\theta}^2d\theta+\left(\gamma k^4-\frac{2\lambda \eta}{\nu}k^4\right)\frac{1}{2P}\int_{-P}^P u_{0,\theta\theta}^2 d\theta.
\end{split}
\edis

Next, omitting the details, it can be shown that by manipulating (\ref{1oi2}),
\bdis\bspl
&\frac{1}{2P}\frac{\lambda k^2}{\nu}\int_{-P}^P u_{0,\theta}^2 d\theta=\frac{2}{5}(3\hat{D}+2\hat{C}d+U(\hat{C}+Ud)),\\
&\frac{1}{2P}\frac{\lambda^2k^4}{\nu^2}\int_{-P}^Pu_{0,\theta\theta}^2=\frac{1}{7}\left(6\hat{D} d +8 \hat{C}^2+U(-6\hat{D}+6\hat{C}d+2U\hat{C}+2U^2d)\right).
\end{split}\edis
As well, since $u_0$ and its derivatives with respect to $\theta$ are periodic over $2P$, then
\bdis
\int_{-P}^P \bar{V}(u_0)d\theta=\int_{-P}^P(\mu k^2u_{0,\theta\theta}+\gamma k^4 u_{0,\theta\theta\theta\theta}+\eta k^2(u_0^2)_{\theta\theta})d\theta=0.
\edis

Thus,
\beq
M_i=-\frac{k\int_{-P}^Pu_0\bar{V}(u_0)d\theta}{2P\partial_{r_i}k \prod_{i\ne j}(r_i-r_j)},\label{Midef}\eeq
where
\beq\bspl
\frac{1}{2P}\int_{-P}^Pu_0\bar{V}(u_0)d\theta=&-\left(\mu k^2+2U\eta k^2\right)\frac{2\nu}{5\lambda k^2}(3\hat{D}+2\hat{C}d+U(\hat{C}+Ud))\\
&+\left(\gamma k^4-\frac{2\lambda \eta}{\nu}k^4\right)\frac{\nu^2}{7\lambda^2 k^4}\left(6\hat{D} d +8 \hat{C}^2+U(-6\hat{D}+6\hat{C}d+2U\hat{C}+2U^2d)\right).
\end{split}\label{diffint}\eeq

\section{Steady Solutions of the Modulation Equations}

The differential equations have been derived that govern the modulation of our leading order cnoidal wave solution (\ref{u0defc5}). We now seek steady solutions to the system (\ref{modsys})-(\ref{Midef}) by setting the wave speed, $c$, to a constant. Here, the workings of \citet{el05} are followed, where steady solutions were found of the fourth order Whitham system for the Kaup-Boussinesq-Burgers equation. This method is now adapted to analyse our third order system.

For $c$ is some constant,
\bdis
\frac{\partial\;}{\partial T}=-\omega\frac{d\;}{d{\Theta}},\quad\frac{\partial\;}{\partial X}=k\frac{d\;}{d{\Theta}}.\edis
Also,
\bdis
\frac{\partial r_i}{\partial T}+Q_i \frac{\partial r_i}{\partial X}=(-k\nu U+k\nu U +k\tilde{Q}_i)\frac{d r_i}{d {\Theta}}=M_i,\edis
where
\bdis
\tilde{Q}_i=Q_i-\nu U=\frac{\nu}{6}\frac{k}{\partial_{r_i}k},\quad M_i=-\frac{k\tilde{M}}{2\partial_{r_i}k \prod_{i\ne j}(r_i-r_j)},
\quad\tilde{M}=2\left(\frac{1}{2P}\int_{-P}^Pu_0\bar{V}(u_0)d\theta\right).\edis
Refer to (\ref{diffint}) to write $\tilde{M}$ in full. Therefore,
\begin{equation}
\frac{d r_i}{d{\Theta}}=\frac{M_i}{k\tilde{Q}_i}=\frac{\hat{M}}{\Pi _{i\ne j}(r_i-r_j)},\quad\hat{M}=-\frac{3}{k\nu}\tilde{M}.\label{rith}\eeq

Now let
\begin{equation}
P(r)=\Pi_{i=1}^3(r-r_i)=r^3-s_1 r^2+s_2 r-s_3,\label{Pdef}\eeq
where
\begin{equation}
s_1=r_1+r_2+r_3,\quad s_2=r_1r_2+r_1r_3+r_2r_3,\quad s_3=r_1r_2r_3.\eeq
Then,
\bdis
\frac{d s_1}{d{\Theta}}=\frac{d\;}{d{\Theta}}(r_1+r_2+r_3)
=\left(\sum_{i=1}^3\frac{1}{\Pi_{i\ne j}(r_i-r_j)}\right)\hat{M}
=0.
\edis
As well,
\bdis
\frac{d s_2}{d{\Theta}}=\frac{d }{d{\Theta}}(r_1r_2+r_1r_3+r_2r_3)
=\left(\sum_{i=1}^3\frac{\sum_j'r_j}{\Pi_{i\ne j}(r_i-r_j)}\right)\hat{M}
=0.
\edis
Hence, $s_1$ and $s_2$ are any real constants. Next,
\begin{equation}
\frac{d s_3}{d{\Theta}}=\frac{d\;}{d{\Theta}}(r_1r_2r_3)
=s_3\left(\sum_{i=1}^3\frac{1}{r_i}\cdot\frac{1}{\Pi_{i\ne j}(r_i-r_j)}\right)\hat{M}
=s_3\cdot\frac{1}{s_3}\cdot\hat{M}
=\hat{M}.
\label{c5s3thmh}\eeq

The three identities
\bdis
\sum_{i=1}^n\frac{1}{\Pi_{i\ne j}(r_i-r_j)}=0,\quad
\sum_{i=1}^n\frac{\sum_j'r_j}{\Pi_{i\ne j}(r_i-r_j)}=0,\quad
\sum_{i=1}^n\frac{1}{r_i}\cdot\frac{1}{\Pi_{i\ne j}(r_i-r_j)}=\frac{(-1)^{n-1}}{\Pi_i^n r_i},\edis
were used to determine $s_{i,\Theta}$, which are from \citet{el05}.

Now, we know that
\bdis\bspl
\frac{\lambda k^2}{\nu}u_{0,\theta}^2
=&-\frac{1}{3}(u_0-\frac{1}{2}(r_1+r_2-r_3))(u_0-\frac{1}{2}(r_1+r_3-r_2))(u_0-\frac{1}{2}(r_2+r_3-r_1)),\\
=&2\hat{D}+2\hat{C}u_0+Uu_0^2-\frac{1}{3}u_0^3.
\end{split}\edis
This becomes, if $u_0=-\tilde{u}_0+s_1/2=-\tilde{u}_0+(r_1+r_2+r_3)/2$ and using (\ref{Pdef}),
\begin{equation}\bspl
\frac{\lambda k^2}{\nu}\tilde{u}_{0,\theta}^2=&\frac{1}{3}(\tilde{u}_0-r_3)(\tilde{u}_0-r_2)(\tilde{u}_0-r_1)\\
=&\frac{1}{3}P(\tilde{u}_0) =\frac{1}{3}(\tilde{u}_0^3-s_1\tilde{u}_0^2+s_2\tilde{u}_0-s_3),\\
=&2\hat{D}+2\hat{C}(-\tilde{u}_0+s_1/2)+U(-\tilde{u}_0+s_1/2)^2-\frac{1}{3}(-\tilde{u}_0+s_1/2)^3.\end{split}\label{uotc5NS}\eeq
Equating like terms of $\tilde{u}_0^n$ in (\ref{uotc5NS}), we find
\begin{equation}
U=\frac{s_1}{6},\quad\hat{C}=-\frac{s_2}{6}+\frac{s_1^2}{24},\quad\hat{D}=-\frac{s_3}{6}+\frac{s_1}{48}(4s_2-s_1^2).\label{somdef}\eeq

The integration constant, $\hat{C}$, can be written as a function of $a$ and $m$ only using (\ref{abddefb}) and (\ref{abddefc}). After some detail, (\ref{cdeqna}) takes the form
\bdis
\hat{C}
=\frac{a^2}{18m^4}(m^4-m^2+1)-\frac{s_1^2}{72}.\edis
Therefore, from (\ref{somdef}),
\bdis
-\frac{s_2}{6}+\frac{s_1^2}{24}=\frac{a^2}{18m^4}(m^4-m^2+1)-\frac{s_1^2}{72}.\edis
Rearranging this,
\begin{equation}
a(m)=\kappa^{1/2}\sqrt{\frac{18m^4}{m^4-m^2+1}},\label{adefc5}\eeq
where
\begin{equation}
\kappa=\frac{s_1^2}{18}-\frac{s_2}{6}>0,\label{c5kapdef}\eeq
which is some constant since $s_1,\;s_2$ are constants. Consequently, (\ref{abddefb}) takes the form
\begin{equation}
d(m)=\frac{s_1}{6}-\frac{\sqrt{\kappa}}{3}\sqrt{\frac{18}{m^4-m^2+1}}\left(2-m^2-3\frac{E(m)}{K(m)}\right).
\label{ddefc5NS}\eeq
Next, from (\ref{abddefc}), (\ref{adefc5}) and (\ref{ddefc5NS}), we can express $\hat{D}$ as a function of $m$ only. With some simplification, (\ref{cdeqnb}) can be written
\bdis
\hat{D}
=\frac{\sqrt{2}}{3}\left(\frac{\kappa}{m^4-m^2+1}\right)^{3/2}(-2m^6+3m^4+3m^2-2)+\frac{s_1s_2}{36}-\frac{11s_1^3}{1296}.\edis
Combining this with our second definition for $\hat{D}$, see (\ref{somdef}),
\bdis
-\frac{s_3}{6}+\frac{s_1}{48}(4s_2-s_1^2)=\frac{\sqrt{2}}{3}\left(\frac{\kappa}{m^4-m^2+1}\right)^{3/2}(-2m^6+3m^4+3m^2-2)+\frac{s_1s_2}{36}-\frac{11s_1^3}{1296}.\edis
This gives us, after some manipulation, an expression for $s_3$ in terms of $m$ only,
\begin{equation}
s_3(m)=-\frac{1}{27}\left(\frac{s_1^2-3s_2}{m^4-m^2+1}\right)^{3/2}(-2m^6+3m^4+3m^2-2)+\frac{s_1s_2}{3}-\frac{2s_1^3}{27},\eeq
and then,
\begin{equation}
\frac{ds_3}{dm}=\frac{(s_1^2-3s_2)^{3/2}m^3(m^2-1)}{(m^4-m^2+1)^{5/2}}.\label{s3m}\eeq

Thus, from (\ref{c5s3thmh}),
\beq\bspl
\frac{dm}{d{\Theta}}=-\frac{3}{k\nu}\left(\frac{ds_3}{dm}\right)^{-1}\tilde{M}
=-\frac{3}{k\nu}\frac{(m^4-m^2+1)^{5/2}}{(s_1^2-3s_2)^{3/2}m^3(m^2-1)}\tilde{M},
\end{split}\label{dmdT}\eeq
where
\beq\bspl
\tilde{M}=&-\left(\mu k^2+2U\eta k^2\right)\frac{4\nu}{5\lambda k^2}(3\hat{D}+2\hat{C}d+U(\hat{C}+Ud))\\
&+\left(\gamma k^4-\frac{2\lambda \eta}{\nu}k^4\right)\frac{2\nu^2}{7\lambda^2 k^4}\left(6\hat{D} d +8 \hat{C}^2+U(-6\hat{D}+6\hat{C}d+2U\hat{C}+2U^2d)\right).\end{split}\label{tilMdef}\eeq
The differential equation (\ref{dmdT}) describes the slow modulations of the cnoidal wave solution (\ref{u0defc5}) as it propagates with constant speed $c$, where $\theta=\Theta/\epsilon$. However, as we seek solutions that satisfy periodic boundaries, one way to achieve this is to set $m_{\Theta}=0$, and therefore, $\tilde{M}=0$. As a result, the modulation term $m$, and therefore, the wave amplitude and period will remain constant over the solution domain. \citet{hat15} used the technique of fixing $m$ to obtain spatially periodic solutions to the periodically forced steady KdV-Burgers equation. The period of these solutions was equal to or an integer multiple of the forcing term's period.

Now, let us express our parameters in terms of $m$, $s_1$ and $\kappa$ such that
\beq
a=\sqrt{\kappa}H_1(m),\quad
d=\frac{s_1}{6}+\sqrt{\kappa}H_3(m),\quad
\hat{C}=\kappa-\frac{s_1^2}{72},\quad
\hat{D}=\kappa^{3/2}H_2(m)-\frac{s_1}{6}\left(\kappa-\frac{s_1^2}{216}\right),\label{adCDdefs}\eeq
where
\bdis\bspl
&{H_1}(m)=\left(\frac{18m^4}{1-m^2+m^4}\right)^{1/2},\\
&H_2(m)=\frac{{H_1}(m)^3}{162 m^6}(-2+3m^2+3m^4-2m^6),\\
&H_3(m)=\frac{H_1(m)}{3 m^2}\left(\frac{3E(m)}{K(m)}+m^2-2\right).
\end{split}
\edis
Substituting the definitions given by (\ref{adCDdefs}) into (\ref{tilMdef}) and simplifying, (\ref{tilMdef}) reduces to
\bdis
\tilde{M}
=-2\left(\mu+\frac{s_1}{3}\eta \right) \frac{2\nu\kappa^{3/2}}{5\lambda}(3H_2(m)+2H_3(m))+2\left(\gamma-\frac{2\lambda \eta}{\nu}\right)\frac{\nu^2\kappa^2}{7\lambda^2}(6H_2(m)H_3(m)+8).
\edis
Therefore, to ensure $m_{\Theta}=0$, then
\beq
-\left(\mu+\frac{s_1}{3}\eta \right) \frac{2\nu\kappa^{3/2}}{5\lambda}(3H_2(m)+2H_3(m))+\left(\gamma-\frac{2\lambda \eta}{\nu}\right)\frac{\nu^2\kappa^2}{7\lambda^2}(6H_2(m)H_3(m)+8)=0.\label{fpcrit}\eeq
If (\ref{fpcrit}) holds then the leading order solution will be periodic since it becomes the cnoidal wave of constant modulus $m$ i.e. the period, amplitude and mean value remain unchanged for all $\theta$.

Given $m$ and $c$ are some constants, then $\omega$ and $k$ will also be constants ($\Theta=k X-\omega T$), where
\beq
k(m)=\sqrt{\frac{\nu a({m})}{12\lambda}}\frac{P}{m K({m})}.
\eeq
This gives, with (\ref{adefc5}),
\beq
\kappa=\frac{s_1^2}{18}-\frac{s_2}{6}=\left(\frac{k K(m)}{P}\right)^4\left(\frac{12\lambda}{\nu}\right)^2\frac{(m^4-m^2+1)}{18}.\label{kapdef}
\eeq

Rearranging (\ref{fpcrit}) and using (\ref{kapdef}), we obtain
\beq
s_1=\frac{180}{7\bar{\rho}(m)}\left(\frac{k}{P}\right)^2\left(\frac{\gamma}{\eta}-\frac{2\lambda}{\nu}\right)-\frac{3\mu}{\eta},\label{s11}\eeq
where
\bdis
\bar{\rho}(m)=\frac{H_1(m)}{(mK(m))^2}\left(\frac{3H_2(m)+2H_3(m)}{3H_2(m)H_3(m)+4}\right).
\edis
This relation determines the speed $c=(\nu s_1)/6$ of our cnoidal wave solutions for some choice of fixed $m$, $(k/P)$ and equation parameters $\lambda,\;\nu,\;\mu,\;\gamma,\;\eta.$ Hence, a family of periodic solutions to (\ref{meqnuns})-(\ref{defVb}) exist in the limit $0<\epsilon\ll 1$ when (\ref{s11}) holds.

\section{Application to the Traffic Flow Model}
The analysis outlined in Sections $3$ and $4$ is now applied to the traffic flow problem discussed in Section $2$. It is necessary to express the leading order solution in terms of the model (\ref{traf_hw}) variables, $j$ and $t$. As well, the equation parameters $\lambda,\nu,\mu,\gamma,\eta$ are now defined using (\ref{par_val}).

Without any loss of generality, we firstly set $k=1$ ($\theta=x$) and then, $c=\omega.$ So, our travelling wave solution in terms of car $j$'s headway at time $t$ is
\beq
\Delta x_j(t)=h+\frac{\epsilon^2}{V^{''}(h)}u_0(j,t) +O(\epsilon^3),\label{hwdef}
\eeq
where, using (\ref{jttoxtb}),
\beq\bspl
u_0(j,t)=&a(m)b(m)+d(m)\\
&\quad +a(m)\operatorname{cn}^2\left(\frac{K(m)}{P}\left(\sqrt{{12(\tau_s-\tau)}}\left(j+t\left(V^{'}(h)+\frac{ s_1}{6}\left(\frac{\tau_s-\tau}{\tau_s}\right)\right)\right)+\theta_0\right);m\right).\end{split}\label{u0jt}
\eeq
The solution parameter $s_1$ is chosen such that (\ref{s11}) is satisfied and $m$ is fixed. The solution parameters $b,\;a$ and $d$ are defined by (\ref{abddefc}), (\ref{adefc5}) and (\ref{ddefc5NS}) respectively. The speed of this wave is
\beq
\textrm{wave speed}=V^{'}(h)+\frac{ s_1}{6}\left(\frac{\tau_s-\tau}{\tau_s}\right).\label{ws_def}\eeq

Note that the solutions highlighted in Section $4$ are periodic in the $\theta$ direction, however, they are not necessarily periodic over the spatial domain $j\in[0,N]$. Since it is periodic boundaries that are of interest here, we set
\bdis
n \hat{T}=N,
\edis
where $n$ is some integer that corresponds to the number of oscillations over $j\in[0,N]$ and $\hat{T}=(2P)/\sqrt{12(\tau_s-\tau)}$ is the solution period in the $j$ direction. Hence,
\beq
P=\frac{N}{2n}\sqrt{12(\tau_s-\tau)}.\label{pdef}\eeq
If $P$ is defined using (\ref{pdef}) then (\ref{pbcon}) holds.

There exist a large number of possible travelling wave solutions that satisfy (\ref{s11}) and (\ref{pdef}). To restrict our analysis, we will seek solutions similar to that found numerically by \citet{zhu08}, such that
\beq
\Delta x_{0,N}(0)=h,\quad \frac{\partial\Delta x_{0,N}(0)}{\partial j}=0,\label{u0ic}\eeq
and therefore,
\bdis
u_0(j=0,N;t=0)=0,\quad \frac{\partial u_{0}}{\partial j}(j=0,N;t=0)=0.\edis

So that (\ref{u0ic}) is met, we set $\theta_0=P$ ($\theta_0$ is arbitrary at leading order) and $ab+d=0$, which gives, using (\ref{adefc5}) and (\ref{ddefc5NS}),
\beq
s_1(m)=-6\left(a(m)b(m)+\frac{a(m)}{3m^2}\left(3\frac{E(m)}{K(m)}+m^2-2\right)\right).\label{s12}\eeq

Now, combining (\ref{s11}) and (\ref{s12}), and after some manipulation, we arrive at
\beq\bspl
\tau
&=\frac{2 n^2}{3N^2V^{'}(h)}\left(\frac{15}{7\bar{\rho}(m)}+6 m^2K(m)^2\left(b+\frac{1}{3m^2}\left(\frac{3E(m)}{K(m)}+m^2-2\right)\right)\right)+\frac{1}{2V^{'}(h)}.
\end{split}\label{tau_def}
\eeq
So, $\tau$ and therefore $\hat{a}=1/\tau$ are functions of $(n/N)$, $m$ and $V^{'}(h)$. Figure \ref{pp_h35} depicts $\hat{a}(m)$ for $h=3.5$ or $4.5$ and $h=2.5$ or $5.5$. Each curve signifies when a cnoidal wave solution satisfying (\ref{u0ic}) occurs for some value of integer $n$ and when $N=100$ (for $100$ cars on the road). Only when $\hat{a}>\hat{a}_s$ will a solution exist ($\epsilon>0$). It is apparent that for a fixed $\hat{a}$, as $n$ becomes larger, the modulation term $m$ decreases.

\begin{figure}
\begin{center}
\includegraphics[width=3.2in]{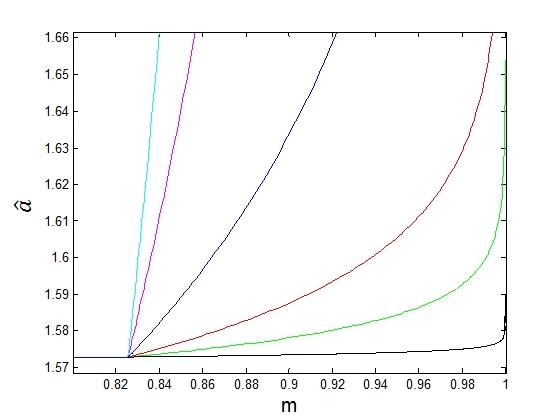}
\includegraphics[width=3.2in]{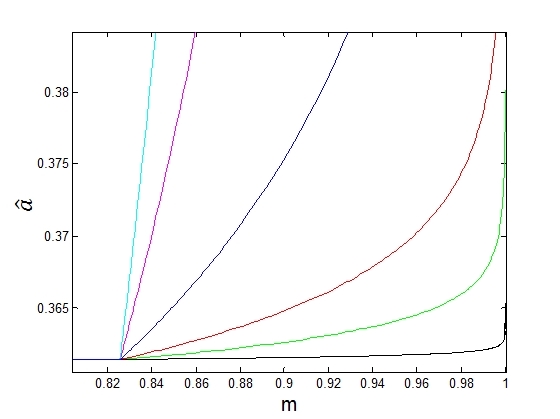}
\caption{The drivers sensitivity, $\hat{a}$, which appears in the OV model (\ref{traf_hw}), is determined by $m$ and $n$. Each curve represents $\hat{a}(m)=1/\tau(m)$ for some choice of $n$, where $\tau$ is given by (\ref{tau_def}), $N=100$, black: $n=1$, green: $n=3$, red: $n=5$, blue: $n=10$, purple: $n=20$, light blue: $n=30$. Left: $h=3.5,\;4.5$, $\hat{a}_s=1.5729$. Right: $h=2.5,\;5.5$, $\hat{a}_s= 0.36141$.
\label{pp_h35}}
\end{center}
\end{figure}

From (\ref{s12}) and (\ref{tau_def}), it is evident that the wave speed (\ref{ws_def}) is determined by also specifying $(n/N)$, $m$ and $V^{'}(h)$. The plot of the wave speed is displayed in Figure \ref{wsfig} for $h=2.5$ or $5.5$ and $h=3.5$ or $4.5$ as a function of the modulus $m$, with each curve corresponding to some choice of $n$ and $N=100$. This figure reveals that the wave speed increases with $m$ for some choice of $n$.

\begin{figure}
\begin{center}
\includegraphics[width=3.2in]{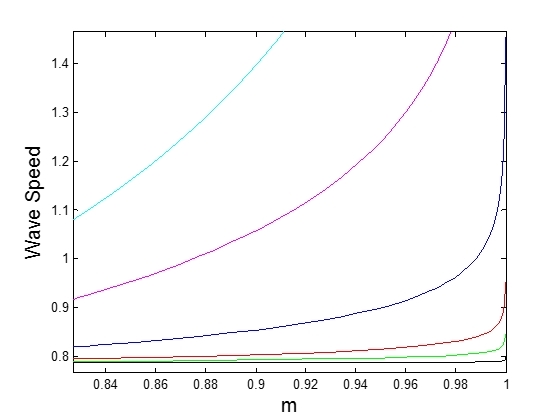}
\includegraphics[width=3.2in]{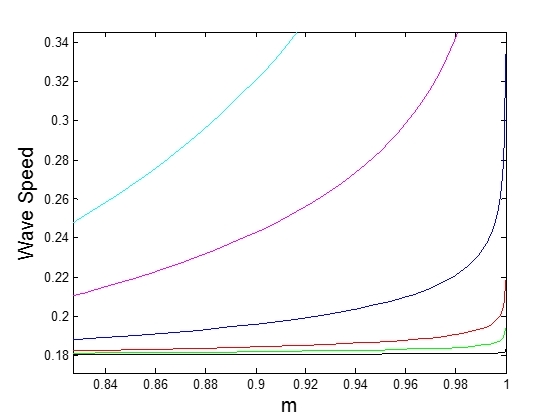}
\caption{The wave speed (\ref{ws_def}) of the cnoidal wave solution (\ref{u0jt}) as a function of the modulus $m$, for some choice of $n$, where $N=100$, black: $n=1$, green: $n=3$, red: $n=5$, blue: $n=10$, purple: $n=20$, light blue: $n=30$. Left: $h=3.5,\;4.5$. Right: $h=2.5,\;5.5$.
\label{wsfig}}
\end{center}
\end{figure}

\section{Results}

The spatially periodic asymptotic solutions, given by (\ref{hwdef})-(\ref{u0jt}), are plotted in Figures \ref{n1}-\ref{n3n5}. So that (\ref{u0ic}) holds, $P$, $s_1$ and $\tau$ are defined using (\ref{pdef}), (\ref{s12}) and (\ref{tau_def}) respectively. As well, the OV model (\ref{traf_hw}) governing the headway is solved numerically with Matlab's ode45, where the periodic boundary conditions (\ref{pbcon}) are imposed. The initial condition used for the simulation is defined by (\ref{hwdef})-(\ref{u0jt}) at $t=0$. Then, we compare the asymptotic solution with the numerical solution at different time intervals.

In Figure \ref{n1}, the solution for $h=3.5,\;\hat{a}=1.59,\;\epsilon=0.10372,\;N=100$ and $n=1$ is shown, where an upward density wave is depicted. The top panel of Figure \ref{n1} compares the asymptotic solution, given by (\ref{hwdef})-(\ref{u0jt}), to the numerical findings for $t\in[0,100]$. The middle panel displays the headway profile for car $j=0,N$, where the solid black curve represents the asymptotic solution and the dotted red curve is the ode45 solution. These plots suggest excellent agreement. Next, the simulation is solved over a large time domain and the result is examined in the bottom panel of Figure \ref{n1} (at around $t=1000$ and $t=10000$), where the headway profiles for car $j=0,N$ is displayed. It is evident that numerically, the density wave persists for a considerably long time, although eventually a slight phase shift and a reduction in amplitude develops. It will eventually disappear. This is expected since all the depicted solutions satisfy the linear stability criteria (\ref{linstab}). Hence, all disturbances will dissolve as $t\rightarrow \infty$ and $\Delta x_j$ will tend to steady state $h$, as discussed in Section $2$. These observations are consistent with the numerical findings of \citet{mur99} for solutions in the metastable zone.

\begin{figure}
\begin{center}
\includegraphics[width=3.2in,height=2.2in]{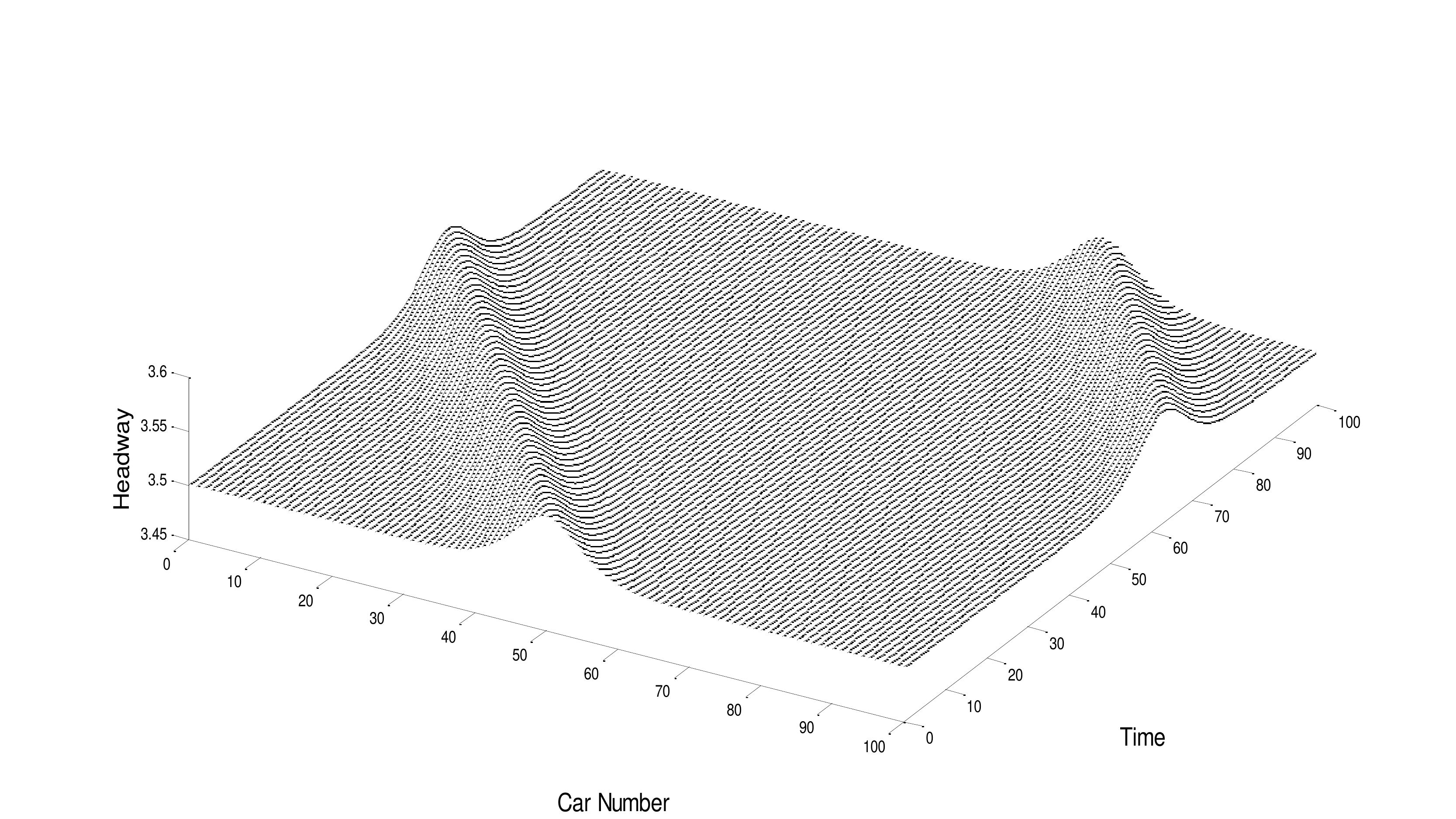}
{\includegraphics[width=3.2in,height=2.2in]{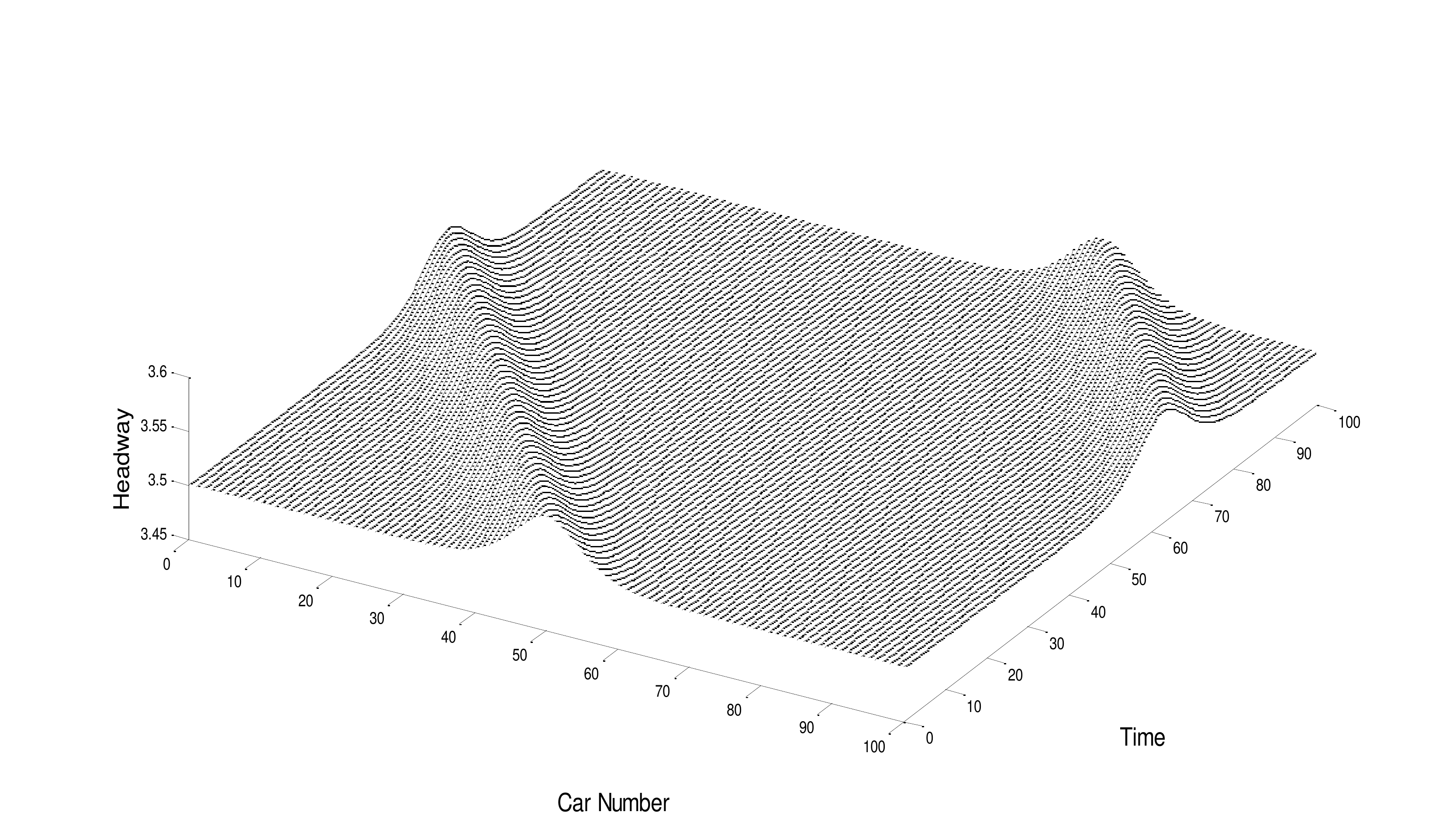}}\\
\includegraphics[width=2.5in]{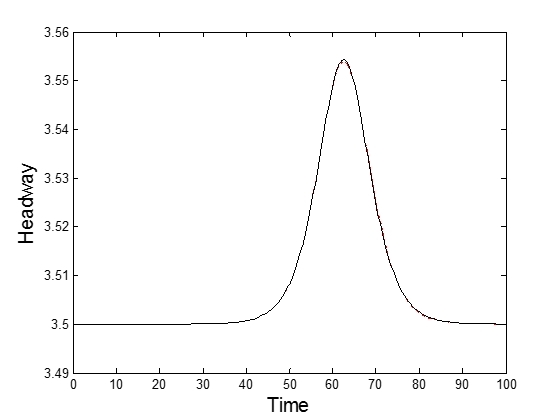}\\
\includegraphics[width=2.5in]{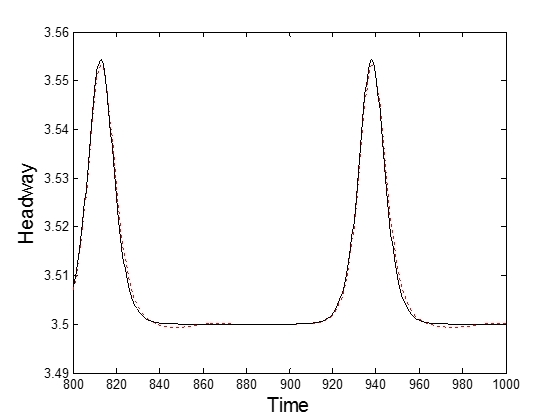}
\includegraphics[width=2.5in]{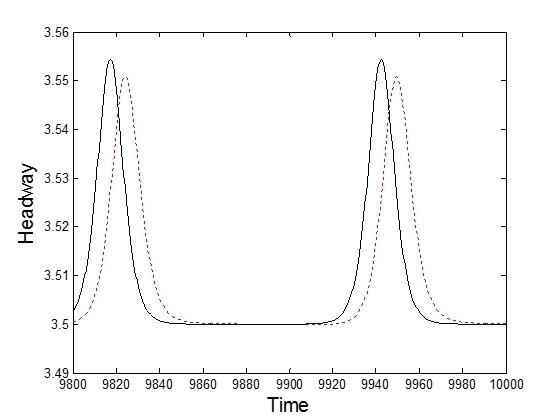}
\caption{Spatially periodic headway solutions for cars $j=0,1,\ldots,100$, with $h=3.5$, $\hat{a}=1.59$, $\epsilon= 0.10372$, $m =0.999998947$, $n =1$, wave speed$=0.79961$. Top left: Asymptotic headway solution given by (\ref{hwdef})-(\ref{u0jt}). Top right: Matlab ode45 headway solution to (\ref{traf_hw}) with the initial condition defined by (\ref{hwdef})-(\ref{u0jt}). Middle: Headway profile for car $j=0,100$, where the asymptotic solution corresponds to the solid black curve and the ode45 solution is represented by the dotted red curve. Bottom: Headway profile for car $j=0,100$ and $t\in[800,1000]$ (left), $t\in[9800,10000]$ (right), where the asymptotic solution corresponds to the solid black curve and the ode45 solution is represented by the dotted red curve.
\label{n1}}
\end{center}
\end{figure}

Next, in Figure \ref{n1c}, the driver's sensitivity is increased to $\hat{a}=1.65$ ($\epsilon=0.21617$), and consequently, the solution is larger in amplitude, with a narrower peak (since $m$ is increased and it is extremely close to $1$). This is of a form similar to a soliton, although, as it is actually a cnoidal wave, the periodic boundary conditions are satisfied. \citet{zhu08} found numerically solutions of this form when they considered periodic boundaries, however they referred to these as solitons. The top panel compares the numerical and asymptotic findings for $t\in[0,100]$. The bottom panel depicts the headway profiles of car $j=0,N$ for $t\in[0,100]$ (left) and $t\in[800,1000]$ (right). The asymptotic and numerical solutions are in good agreement. However, as $\epsilon$ is larger, the excellent match observed for $\hat{a}=1.59$ is not achieved here. Moreover, as a result of increasing $\epsilon$, the solution will disappear significantly faster since a notable phase shift and amplitude reduction appears around $t=900$. This behaviour can be attributed to a greater wave speed.

\begin{figure}
\begin{center}
\includegraphics[width=3.2in,height=2.1in]{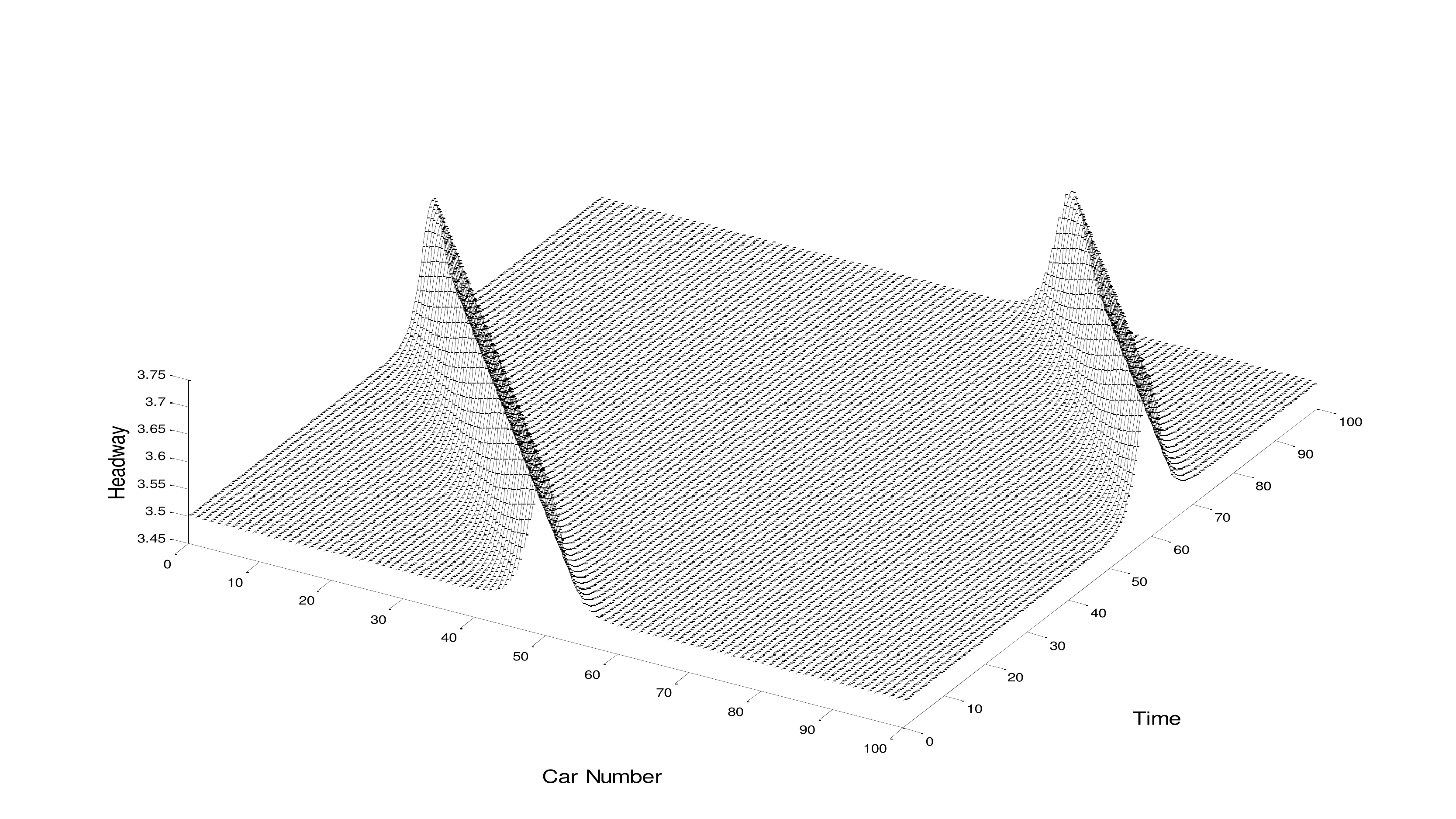}
\includegraphics[width=3.2in,height=2.1in]{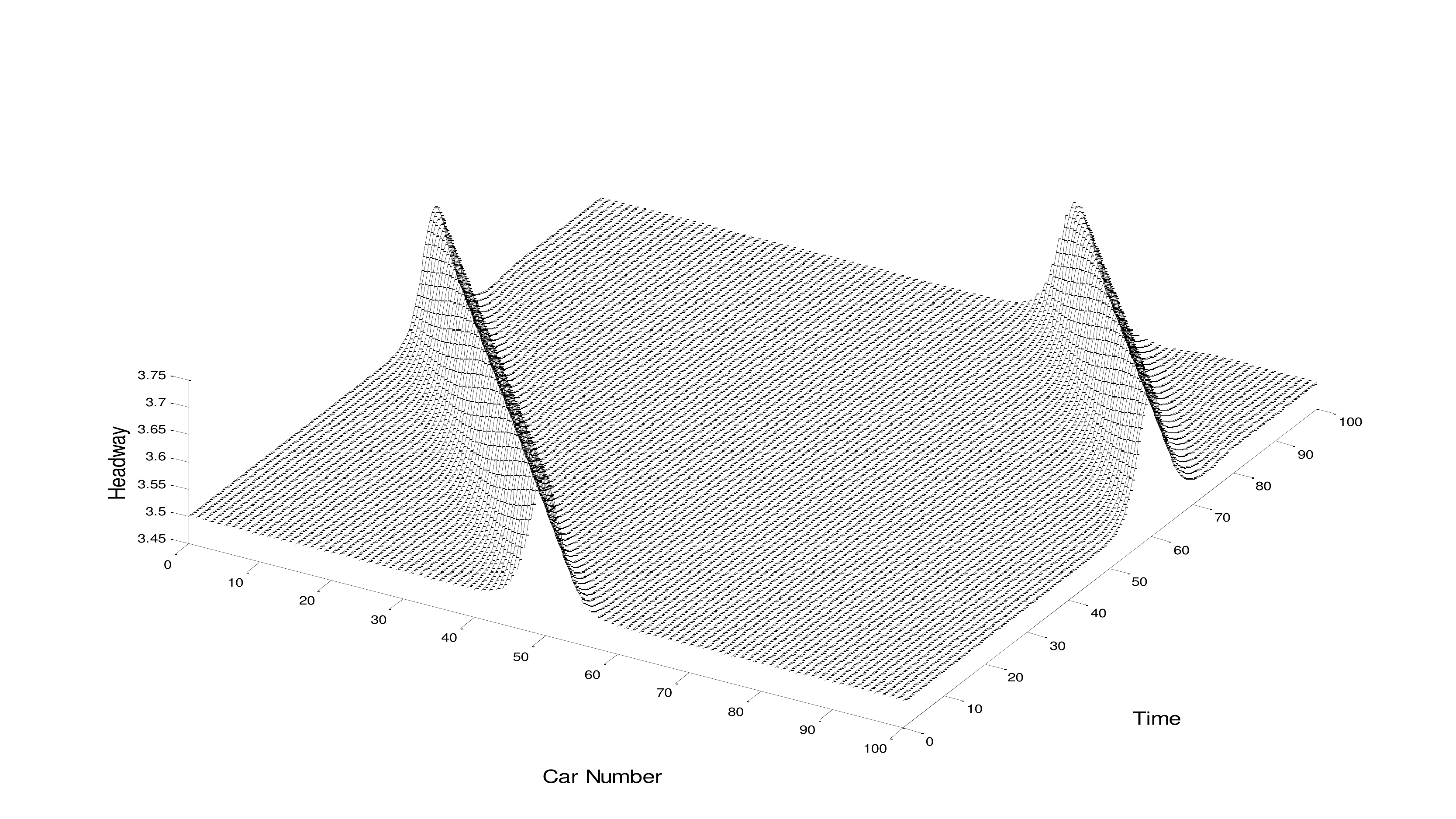}\\
\includegraphics[width=2.5in]{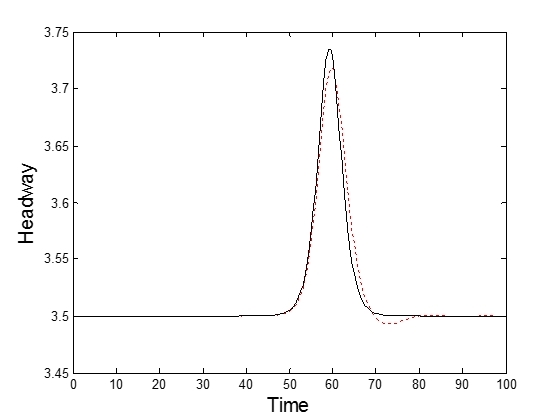}
\includegraphics[width=2.5in]{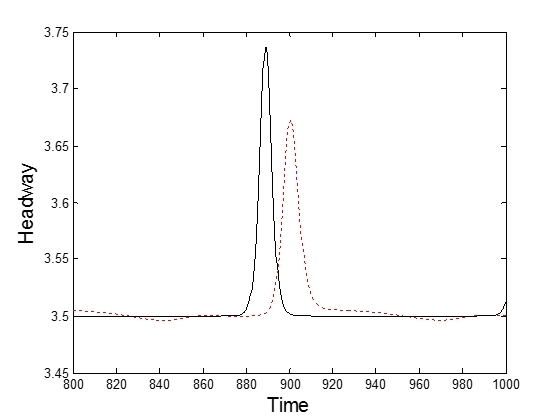}
\caption{Spatially periodic headway solutions for cars $j=0,1,\ldots,100$, with $h=3.5$, $\hat{a}=1.65$, $\epsilon= 0.21617$, $m =0.999999999999963$, $n =1$, wave speed$=0.84357$. Top left: Asymptotic headway solution given by (\ref{hwdef})-(\ref{u0jt}). Top right: Matlab ode45 headway solution to (\ref{traf_hw}) with the initial condition defined by (\ref{hwdef})-(\ref{u0jt}). Bottom: Headway profile for car $j=0,100$ and $t\in[0,100]$ (left), $t\in[800,1000]$ (right), where the asymptotic solution corresponds to the solid black curve and the ode45 solution is represented by the dotted red curve.
\label{n1c}}
\end{center}
\end{figure}

In Figure \ref{n1e}, again solutions for $n=1,\;N=100$ are examined except now $h=4.5$. The same choice of $m$ and wave speed used for $h=3.5$ will apply here. The top panel relates to $\hat{a}=1.59$ ($\epsilon=0.10372$) and the bottom to $\hat{a}=1.65$ ($\epsilon=0.21617$). Asymptotic solutions for $t\in[0,100]$ are shown on the left, and on the right, the asymptotic solutions for car $j=0,N$ are compared to the numerical findings. The agreement between these two solutions and the long-time behaviour are consistent with the discussions for $h=3.5$. However, now the solution is a downward form density wave since $h>h_c=4$. \citet{zhu08} made similar observations, explaining that if $h<h_c$, vehicle $j$ slows down to avoid crashing into vehicle $j+1$. As a result, an upward form density wave occurs. Whereas, if $h>h_c$, vehicle $j$ speeds up to the maximal velocity and consequently, a downward form density wave emerges. Given the traffic density is defined as the inverse of the headway, then the upward and downward waves represent clusters of faster and slower moving vehicles respectively.

\begin{figure}
\begin{center}
\includegraphics[width=3.2in,height=2.1in]{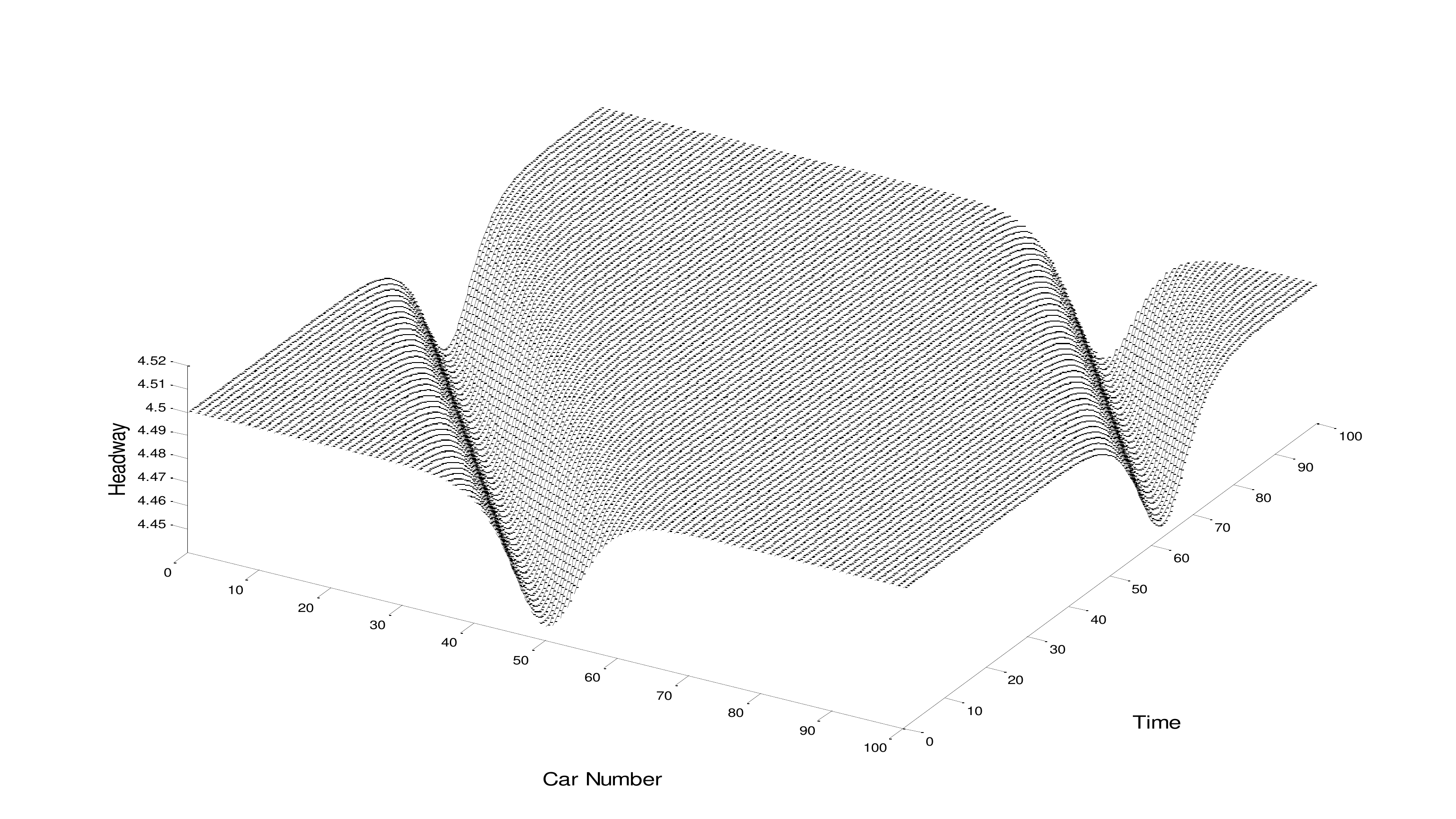}
\includegraphics[width=2.5in]{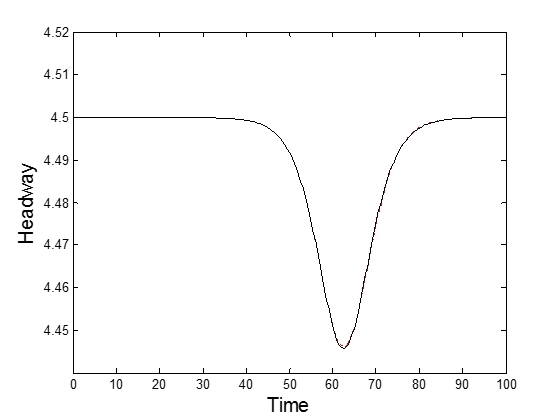}\\
\includegraphics[width=3.2in,height=2.1in]{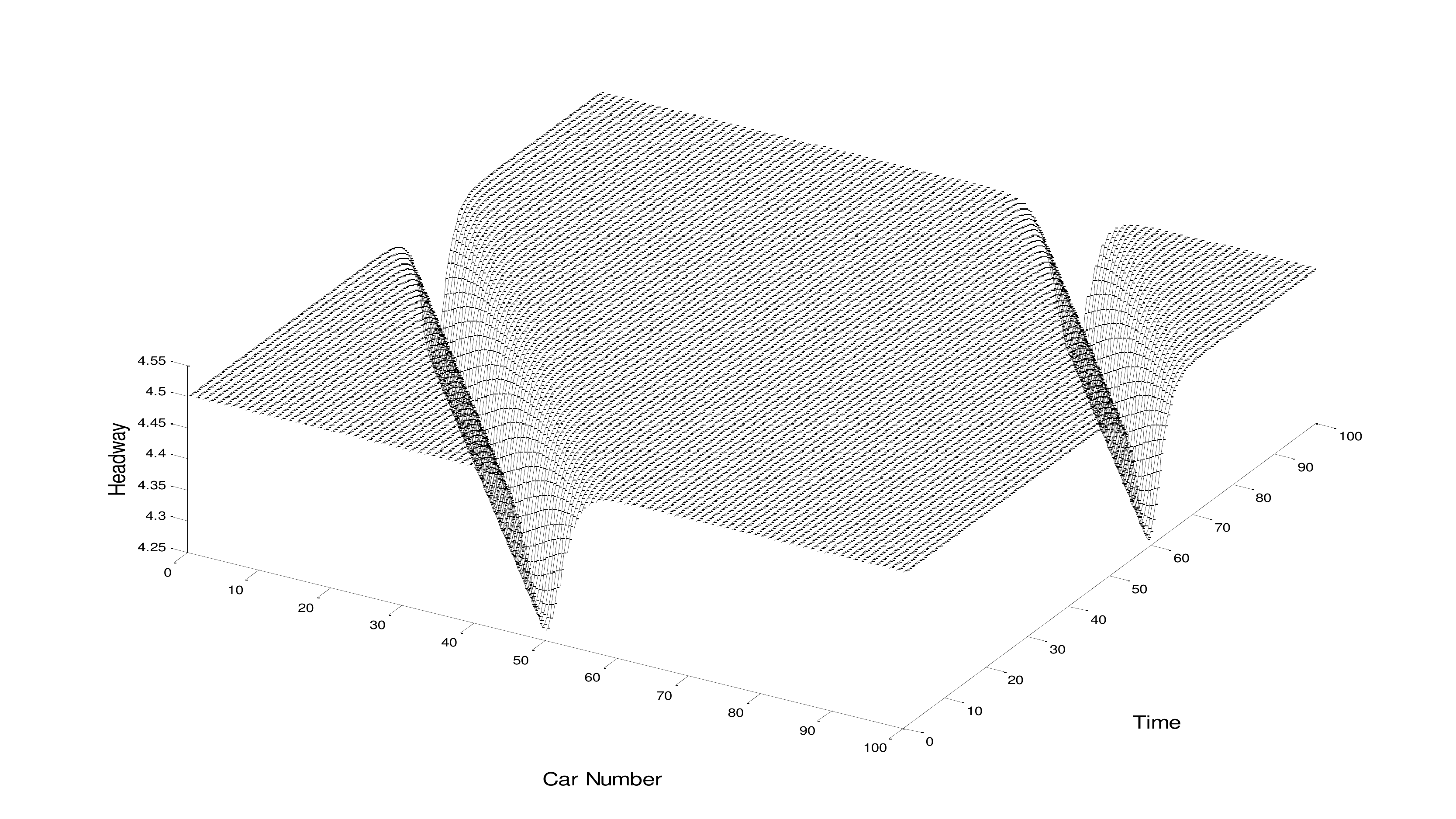}
\includegraphics[width=2.5in]{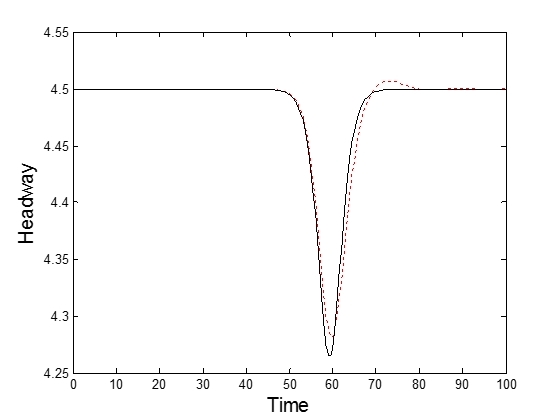}
\caption{Spatially periodic headway solutions for cars $j=0,1,\ldots,100$, with $h=4.5$, $n =1$. Top: $\hat{a}=1.59$, $\epsilon=0.10372$, $m = 0.999998947$, wave speed$=0.79961$. Bottom: $\hat{a}=1.65$, $\epsilon=0.21617$, $m = 0.999999999999963$, wave speed$=0.84357$. Top/Bottom left: Asymptotic headway solution given by (\ref{hwdef})-(\ref{u0jt}). Top/Bottom right: Headway profile for car $j=0,100$, where the asymptotic solution corresponds to the solid black curve and the ode45 solution to (\ref{traf_hw}) is represented by the dotted red curve.
\label{n1e}}
\end{center}
\end{figure}

Now, we investigate the solutions with $n=2$, $h=3.5,\;N=100$ so that two density waves propagate, see Figure \ref{n2}. The top panel corresponds to $\hat{a}=1.59$ ($\epsilon=0.10372$), and the bottom to $\hat{a}=1.65$ ($\epsilon=0.21617$). As exhibited when $n=1$, upward and downward waves form when $h<h_c$ and $h>h_c$ respectively, although solutions with $h>h_c$ are not depicted here. On the left, the asymptotic solution over the domain $t\in[0,100]$ is displayed and on the right, a comparison between the numerical and asymptotic findings for car $j=0,N$ is given. There is good agreement between the two results, especially when $\epsilon=0.10372$, which is expected. Once again, these travelling wave solutions will eventually dissolve at very large $t$, with the $\hat{a}=1.59$ result persisting for a far greater time due to a smaller wave speed. Furthermore, the amplitude is notably larger and narrower for $\epsilon=0.21617$, which is a result of increasing $m$.

\begin{figure}
\begin{center}
\includegraphics[width=3.2in,height=2.1in]{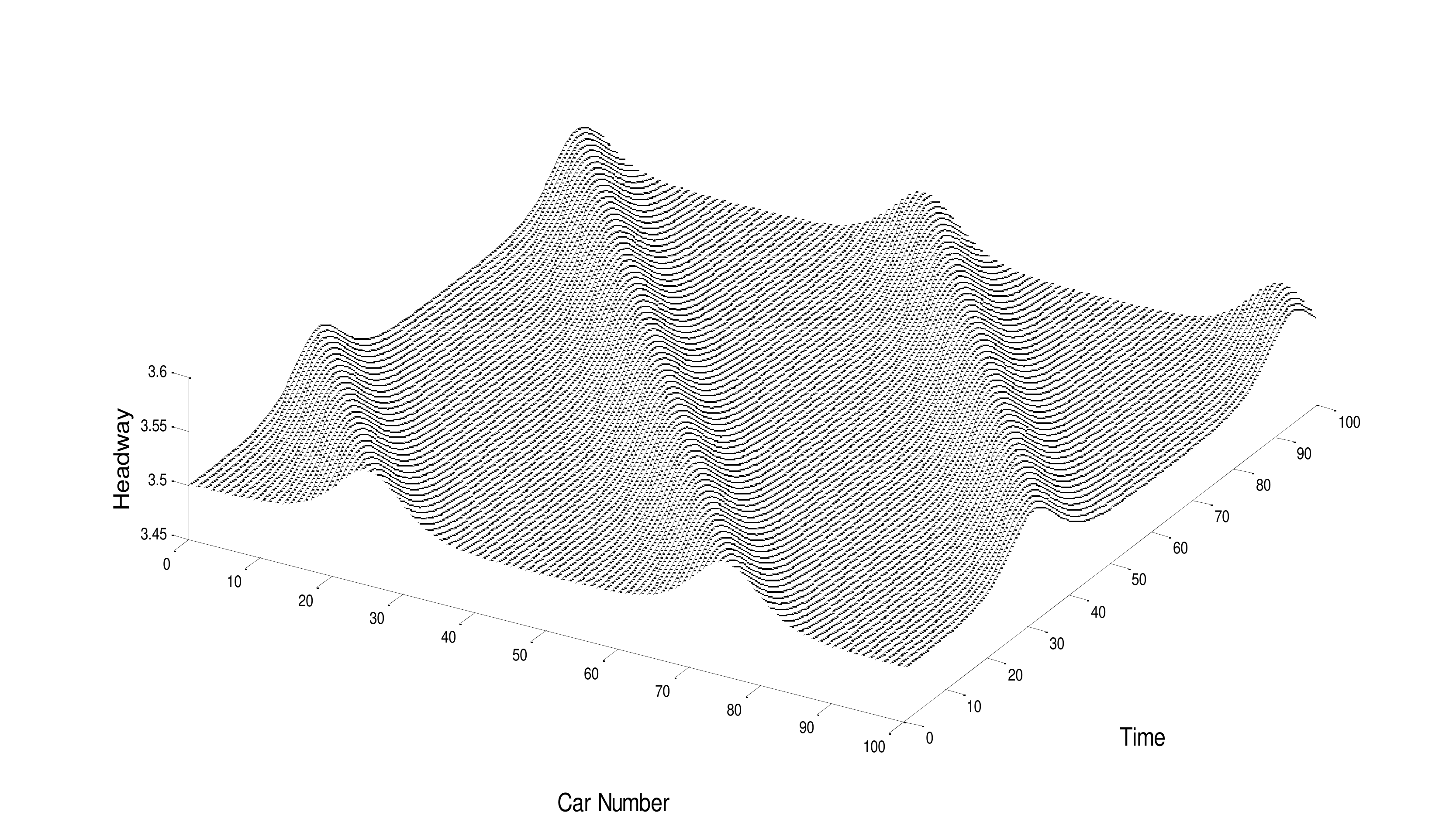}
\includegraphics[width=2.5in]{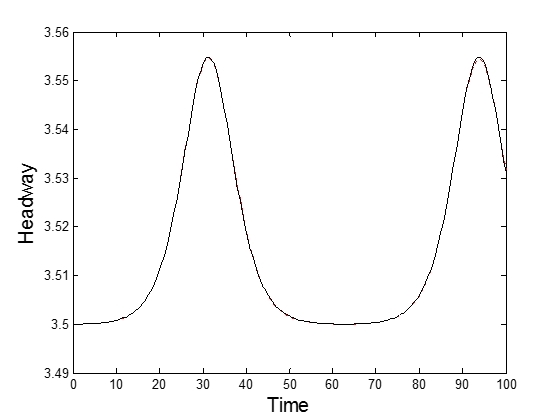}\\
\includegraphics[width=3.2in,height=2.1in]{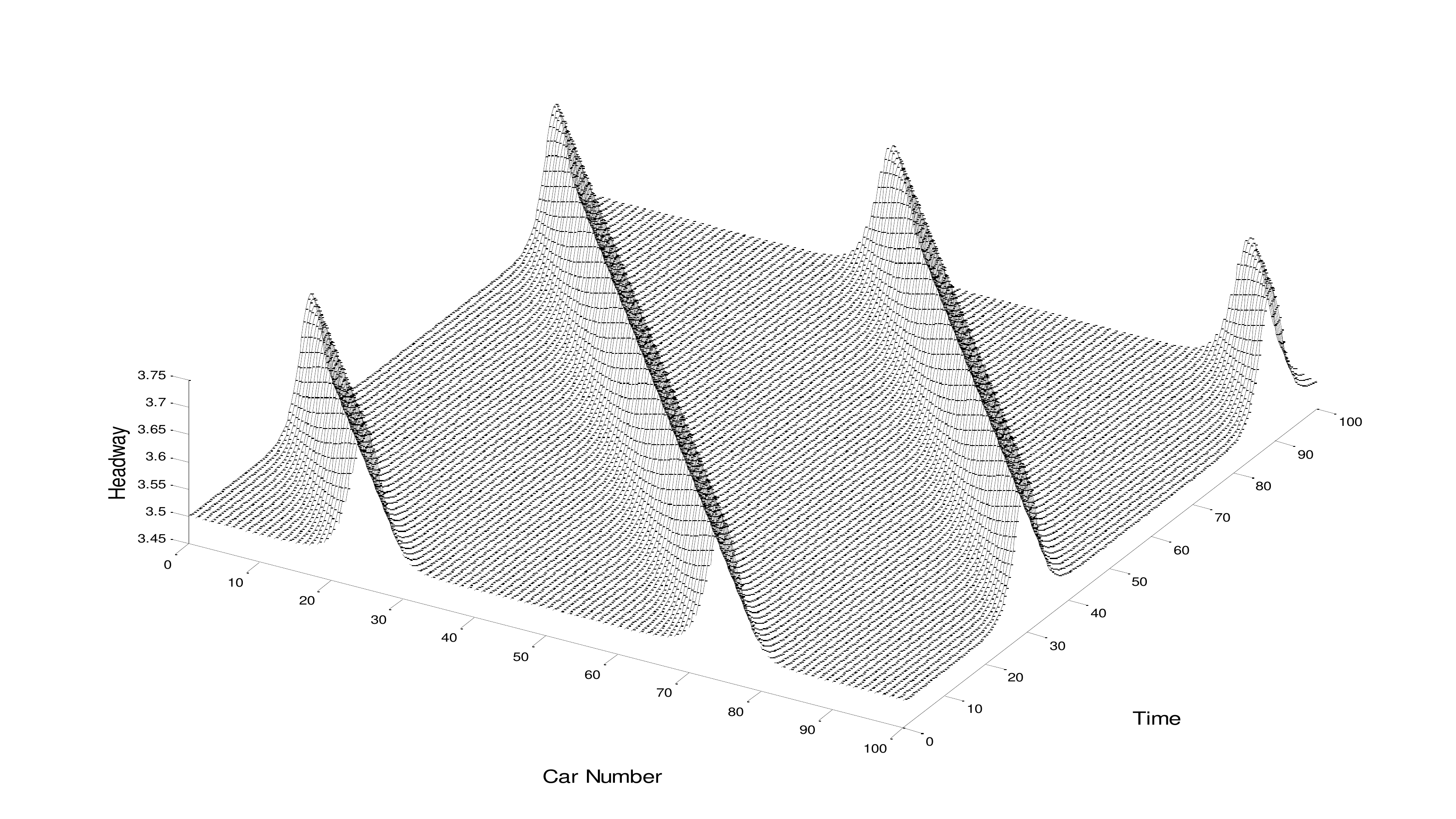}
\includegraphics[width=2.5in]{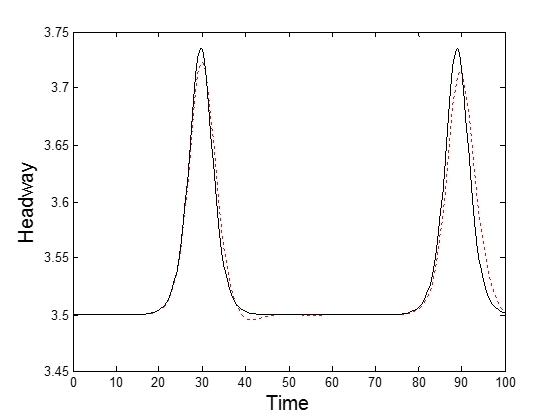}
\caption{Spatially periodic headway solutions for cars $j=0,1,\ldots,100$, with $h=3.5$, $n=2$. Top: $\hat{a}=1.59$, $\epsilon=0.10372$, $m =0.99724797$, wave speed$=0.79967$. Bottom: $\hat{a}=1.65$, $\epsilon=0.21617$, $m = 0.99999946$, wave speed$=0.84362$. Top/Bottom left: Asymptotic headway solution given by (\ref{hwdef})-(\ref{u0jt}). Top/Bottom right: Headway profile for car $j=0,100$, where the asymptotic solution corresponds to the solid black curve and the ode45 solution to (\ref{traf_hw}) is represented by the dotted red curve.
\label{n2}}
\end{center}
\end{figure}
Such studies as \citet{yu10} highlighted numerical density waves with multiple oscillations over the same spatial domain. To obtain this solution type, $n$ can continue to be increased. Analysing these solutions, it becomes apparent that as $n$ grows, for some fixed $\hat{a}$, the amplitude and wave speed slightly increases. As well, these solutions are found to have a similar long-time behaviour as that demonstrated for $n=1$. As an example, Figure \ref{n3n5} depicts an asymptotic solution with three oscillations and where $\hat{a}=1.59,\;N=100,\;h=3.5$. 
\begin{figure}
\begin{center}
\includegraphics[width=3.2in,height=2.1in]{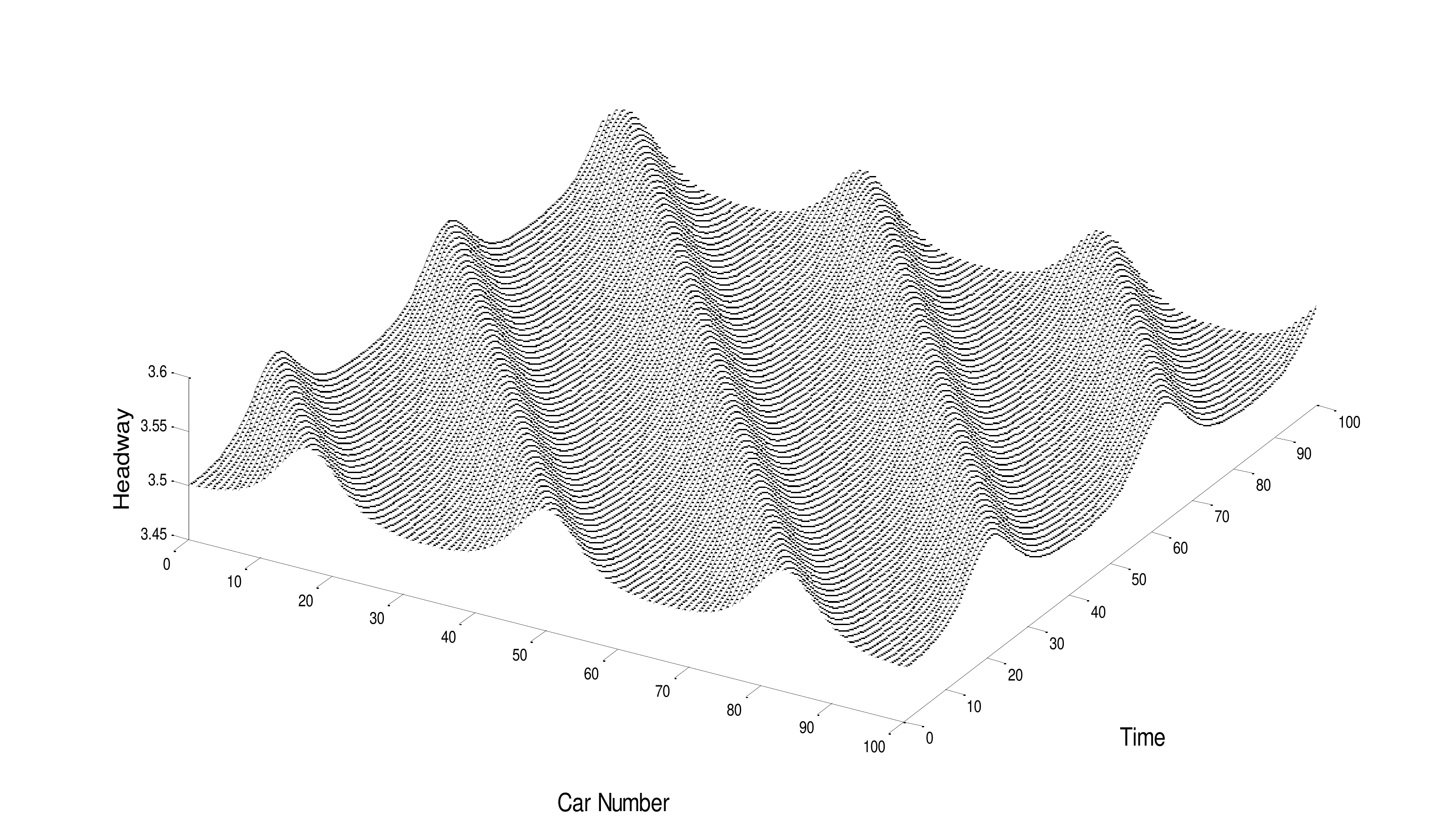}
\includegraphics[width=2.5in]{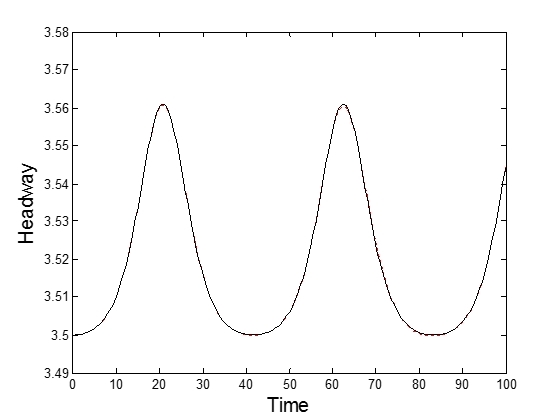}
\caption{Spatially periodic headway solutions for cars $j=0,1,\ldots,100$, with $h=3.5$, $\hat{a}=1.59$, $\epsilon=0.10372$, $n=3$, $m =0.9728972$, wave speed$=0.80039$. Left: Asymptotic headway solution given by (\ref{hwdef})-(\ref{u0jt}). Right: Headway profile for car $j=0,100$, where the asymptotic solution corresponds to the solid black curve and the ode45 solution to (\ref{traf_hw}) is represented by the dotted red curve.
\label{n3n5}}
\end{center}
\end{figure}

\section{Conclusion}
The evolution of traffic behaviour was determined using the OV model (\ref{traf_hw}). This model has been examined numerically and analytically by various previous studies, with the stability regimes well-outlined (refer to \citet{ge05}). Here, we concentrated on the metastable zone that corresponded to the onset of traffic jams and where (\ref{traf_hw}) reduced to the perturbed KdV equation (\ref{meqnuns}). A multi-scale perturbation analysis was then applied to (\ref{meqnuns}). As a result, at leading order, the cnoidal wave solution was obtained and at the next order, the Whitham system was derived, which was altered due to the perturbation terms of (\ref{meqnuns}). Next, steady travelling wave solutions were sought by ensuring the wave speed remained constant. Consequently, the three Whitham equations were transformed into a single equation for the slow variation of the modulation term, $m$, defined by (\ref{dmdT}). Setting (\ref{dmdT}) to zero so that $m$ was fixed over the solution domain then led to the identification of periodic cnoidal wave solutions. This analysis was next applied to the traffic flow problem by defining the leading order solution in terms of car $j=1,\ldots,N$ and imposing the periodic boundary conditions along $j\in[0,N]$. Thus, a family of travelling wave solutions were highlighted, where the choice of $m$ and the number of oscillations over the spatial domain were shown to determine the driver's sensitivity, $\hat{a}$, and the wave speed. Lastly, comparisons between the numerical solutions of the OV model and the asymptotic headway solutions were performed. Overall, a good agreement between the two solutions was observed. Although, the numerical density waves dissolved after a considerable length of time, which is consistent with the linear stability analysis. This paper has provided an extension of other workings that only considered traffic soliton solutions (for example see \citet{mur99} and \citet{zhu08}). Instead here modulation theory was used to establish the existence of cnoidal waves in the traffic model.
\bibliographystyle{plainnat2}
\bibliography{bibliog}

\end{document}